\documentclass[a4paper,12pt]{article} 

\usepackage[fleqn]{amsmath}
\usepackage{amscd,amsthm,amssymb}
\usepackage{color}
\usepackage[all]{xypic}
\usepackage[english]{babel}
\usepackage[latin5]{inputenc}
\makeatletter
\newcommand{\specificthanks}[1]{\@fnsymbol{#1}}

\usepackage{setspace}
\newtheorem{teo}{Theorem}[section]
\newtheorem{coro}[teo]{Corollary}
\newtheorem{lema}[teo]{Lemma}
\newtheorem{propo}[teo]{Proposition}

\theoremstyle{definition}
\newtheorem{defi}[teo]{Definition}
\newtheorem{ejem}[teo]{Example}

\theoremstyle{remark}

\newtheorem{obser}[teo]{Remark}

\numberwithin{equation}{subsection}

\def\C{\mathcal{C}}

\def\R{\mathcal{R}}
\def\RR{\mathbb{R}}

\def\d{\mathrm{d}}

\def\NN{\mathbb{N}}
\def\TT{\mathbb{T}}

\def\qed{\hfill $\square$}
\def\demo{\proof}
\def\proof{\noindent \textit{Proof: }}

\def\Conn{\mathcal{C}}
\def\Diff{\mathrm{Diff}_x}
\def\Gl{{\rm Gl}_n}

\begin{document}

\title{\textsc{On moduli spaces for finite-order jets of linear connections}}

\author{Gordillo, A.\thanks{Universidad de Extremadura (SPAIN)\newline  
Both authors have been partially supported by Junta de Extremadura and FEDER funds.} \and Navarro, J.\textsuperscript{\specificthanks{1}\,,}\thanks{Corresponding author. Email address: \texttt{navarrogarmendia@unex.es}\newline Department of Mathematics, Universidad de Extremadura, Avda. Elvas s/n, 06006, Badajoz, Spain.}}

\date{\today}

\maketitle

\begin{abstract}
We describe the ringed--space structure of moduli spaces of jets of linear connections (at a point) as orbit spaces of certain linear representations of the general linear group.

Then, we use this fact to prove that the only (scalar) differential invariants associated to linear connections are constant functions, as well as to recover various expressions appearing in the literature regarding the Poincar\'{e} series of these moduli spaces.

\bigskip

\noindent \emph{Key words and phrases:} Moduli spaces, jets of linear connections,
differential invariants.

\medskip

\noindent MSC: 58D27, 53A55
\end{abstract}

\tableofcontents

\section*{Introduction} 

The aim of this paper is to study the classification of finite-order jets of linear connections at a point. This problem, as well as similar ones regarding local classifications of other geometric structures, have already been widely discussed in the literature (\cite{Teresa}, \cite{Arnold}, \cite{Dubrovskiy}, \cite{DubrovskiyII}, \cite{GNS}, \cite{Jaime} or \cite{Shmelev}).

A classical approach to these type of questions consists on developing ``normal forms'' for the geometric objects under consideration; that is to say, trying to find suitable coordinate charts where the expression of the objects is particularly simple (e.g., \cite{Teresa},\cite{Opozda}). 

Another point of view tries to consruct moduli spaces; i.e., to determine the structure of the orbit space for the action of the Lie pseudo-group  of ``changes of coordinates'' on the space of objects to classify. Usually, this goal is exceedingly difficult, so that, to tackle it, one restricts his attention to infinitesimal neighbourhoods, and hence to jets of the objects under study (see the programme outlined in \cite{Arnold}, Sect. 1, and the development carried out in \cite{Dubrovskiy}, \cite{GNS} or \cite{Jaime}).

In this paper we adopt this latest approach, and hence study the structure of the quotient $\,J^r_x \mathcal{C} / \Diff $, where $\,J^r_x \mathcal{C}\,$ denotes the space of $r$-jets of linear connections at a point, and $\,\Diff\,$ stands for the group of germs of diffeomorphisms leaving the point fixed. 

Our main result, Theorem \ref{elmoduli}, establishes an isomorphism of ringed spaces between this quotient and the orbit space of a linear representation of the general linear group $\, \Gl$. Such a representation is worked out using the so-called {\it normal tensors} associated to linear connections. These tensors were already used in the early developments of Riemannian geometry (\cite{Thomas}) as well as in the theory of natural operations in Riemannian geometry (\cite{Epstein}, \cite{Stredder}). More recently, it was realized that they are particularly well-behaved in order to construct moduli spaces and they have been used to study the classification of certain $G-$structures, including jets of Riemannian metrics (\cite{GNS}, \cite{Jaime}). This paper applies these techniques to the case of linear connections, and obtains the aforementioned result as a corollay of an orbit-reduction-type statement (Theorem \ref{LemaFundamental}).  


In the last section, we use our main result to prove the absence of non-trivial differential invariants associated to linear connections (Proposition \ref{finitud}). Finally, we also make some comments on how to apply Theorem \ref{elmoduli} to recover certain formulae already appearing in the literature (\cite{Dubrovskiy}), where they were obtained through lengthy computations.

\section{Preliminaries}

\subsection*{Ringed spaces}

Apart from trivial cases, the moduli space to be studied throughout this paper is never a smooth manifold. Nevertheless, it can be endowed with certain geometric structure.

To be precise, it will be considered as a ringed space in the following sense:

\begin{defi}
Let $\,X\,$ be a topological space. 
A \textbf{sheaf of continuous
functions} $\, \mathcal{O}_X\,$ on $\,X\,$ is a sub-sheaf of the sheaf $\, \mathcal{C}\,$ 
of real-valued, continuous functions on $\, X$.
\end{defi}

In other words, a sheaf of continuous functions 
on $\,X\,$ is a map $\,\mathcal{O}_X\,$ which assigns a
subalgebra $\,\mathcal{O}_X(U)\subseteq
\mathcal{C}(U,\mathbb{R})\,$ to every open subset $\,U\subseteq
X\,$, with the following condition:

For every open subset $\,U\subseteq X\,$, every open cover
$\,U=\bigcup U_i\,$ and every function $\,f:U\rightarrow
\mathbb{R}\,$, it is verified
$$
f\in \mathcal{O}_X(U)\,\, \Longleftrightarrow \,\, f|_{U_i}\in
\mathcal{O}_X(U_i)\,,\,\,\,\forall\, i\,.
$$



\smallskip

\begin{defi}
We will call \textbf{ringed space} the pair
$\,(X,\mathcal{O}_X)\,$ formed by a topological space $\,X\,$ and
a sheaf of continuous functions $\,\mathcal{O}_X\,$ on $\,X\,$.

Given two ringed spaces $\,X\,$ and $\,Y\,$, a \textbf{morphism of
ringed spaces} $\,\varphi:X\rightarrow Y\,$ is a continuous map
such that, for every open subset $\,V\subseteq Y\,$, the following
condition is held:
$$
f\in \mathcal{O}_Y(V)\,\, \Longrightarrow \,\, f\circ\varphi \in
\mathcal{O}_X(\varphi^{-1}(V))\,.
$$
A morphism of ringed spaces $\,\varphi:X\rightarrow Y\,$ is said
to be an \textbf{isomorphism} if it has an inverse morphism, that
is, there exists a morphism of ringed spaces $\,\phi:Y\rightarrow
X\,$ verifying $\,\varphi \circ \phi=\mbox{Id}_Y\,$, $\,\phi \circ
\varphi=\mbox{Id}_X\,$.
\end{defi}

\medskip





\begin{ejem} \textbf{(Smooth manifolds)} The space
$\,\mathbb{R}^n\,$, endowed with the sheaf
$\,\mathcal{C}_{\mathbb{R}^n}^{\infty}\,$ of smooth functions, is
an example of ringed space. An $\,n-$smooth manifold is precisely
a ringed space in which every point has an open neighbourhood
isomorphic to
$\,(\mathbb{R}^n,\mathcal{C}_{\mathbb{R}^n}^{\infty})\,$. Smooth
maps between smooth manifolds are nothing but morphisms of ringed
spaces.
\end{ejem}

\smallskip

\begin{ejem}\label{ExQuotient}\textbf{(Quotients by the action of a Lie group)}
Let $\,G\times X\rightarrow X\,$ be a smooth action of a Lie group
$\,G\,$ on a smooth manifold $\,X\,$, and let $\,\pi:X\rightarrow
X/G\,$ be the canonical quotient map.

We will consider on the quotient topological space $\,X/G\,$ the
following sheaf $\,\mathcal{C}_{X/G}^{\infty}\,$ of
``differentiable'' functions:

For every open subset $\,V\subseteq X/G\,$,
$\,\mathcal{C}_{X/G}^{\infty}(V)\,$ is defined to be
$$
\mathcal{C}_{X/G}^{\infty}(V):=\{f:V\longrightarrow \mathbb{R}: f
\circ \pi \in \mathcal{C}^{\infty}(\pi^{-1}(V))\}\,.
$$
Note that there exists a canonical $\,\mathbb{R}-$algebra
isomorphism:
$$
\begin{CD}
\mathcal{C}_{X/G}^{\infty}(V) @=
\mathcal{C}^{\infty}(\pi^{-1}(V))^G\\
f & \longmapsto & f\circ \pi\,.
\end{CD}
$$

The pair $\,(X/G,\mathcal{C}_{X/G}^{\infty})\,$ is an example of
ringed space, which we will call \textbf{quotient ringed space} of
the action of $\,G\,$ on $\,X\,$.

As it would be expected, this space verifies the \textbf{universal
quotient property}: Every morphism of ringed spaces
$\,\varphi:X\rightarrow Y\,$, which is constant on every orbit of
the action of $\,G\,$ on $\,X\,$, factors uniquely through the
quotient map $\,\pi:X\rightarrow X/G\,$, that is, there exists a
unique morphism of ringed spaces $\,\tilde{\varphi}:X/G\rightarrow
Y\,$ verifying $\,\varphi=\tilde{\varphi}\circ \pi\,$.
\end{ejem}


\medskip

\subsection*{Invariant theory of the general linear group}

Let $\,V\,$ be an $\,\mathbb{R}$-vector space of finite dimension $\,n$, and let $\,\Gl \,$ be the Lie group of its $\mathbb{R}$-linear automorphisms.

The Main Theorem of the invariant theory for the general linear group (e. g., \cite{KMS}, Sect. 24) states:

\begin{teo}\label{MainTheoremGl}
Let $\,\mathrm{Hom}_{\Gl}\left( V^* \otimes \stackrel{p}{\ldots} \otimes V^* \otimes V \otimes \stackrel{q}{\ldots} \otimes V \, , \, \RR \right) \,$ denote the vector space of $\,\Gl$-invariant, linear maps:
$$ V^* \otimes \stackrel{p}{\ldots} \otimes V^* \otimes V \otimes \stackrel{q}{\ldots} \otimes V \, \longrightarrow \, \RR \ . $$

It holds:
\begin{itemize}
\item If $\,p = q$, then it is spanned by total contractions:
$$ \phi_\sigma (\omega_1 \otimes \ldots \otimes e_p ) \, := \, \omega_1 (e_{\sigma (1)}) \cdot \ldots \cdot \omega_p (e_{\sigma (p)}) \quad , \quad \sigma \in S_p \ . $$

\item If $\,p \neq q$, then that vector space is zero.
\end{itemize}
\end{teo}

\smallskip

To compute invariant functions on a subspace of tensors, recall that $\,\Gl\,$ is linearly semisimple, and hence the following holds:

\begin{propo}\label{Semisimple}
Let $\,E\,$ and $\,F\,$ be linear representations of $\,\Gl$, and let $\,E'\subset E\,$ be a sub-representation. 
Any equivariant linear map $\,E' \to F\,$ is the restriction of an equivariant linear map $\,E \to F$.
\end{propo}

\smallskip

Finally, let us mention that we will be interested in computing {\it smooth}, invariant functions on a linear representation of $\,\Gl$. To this end, the following theorem, which is a particular case of a general result due to Luna (\cite{Luna}), assures that we can always find a system of generators made of polynomial, invariant maps.

\begin{teo}\label{TeoremaLuna}
Let $\, E\, $ be a linear representation of $\,\Gl$, and let $\,\mathsf{A}^{\Gl}\,$ denote the finitely generated algebra of polynomial, $\,\Gl$-invariant functions $\, E \to \RR$. Let $\, p_1, \ldots, p_k\,$ be a system of generators of $\, \mathsf{A}^{\Gl}$ and let us write $\,p = (p_1 , \ldots , p_k ) \colon E \to \RR^k$.

Then,
$$ \mathcal{C}^\infty (E)^{\Gl} \, = \, p^* \mathcal{C}^\infty (\RR^k) \ . $$
\end{teo}

\smallskip

\section{Moduli spaces of jets of 
linear connections}\label{invariantesdiferenciales}

In the remainder of the paper, $\,X\,$ will always be an
$\,n-$dimensional smooth manifold.

Let $\,\Conn \to X\,$ be the bundle of linear connections over $X$, 
and let $\, \widetilde{\Conn} \to X\,$ be the bundle of symmetric, linear connections.

Let us denote by $\,J^r \C \rightarrow X\,$ the fiber bundle of
$\,r-$jets of linear connections on $\,X$. Its fiber over a point
$\,x_0\in X\,$ will be denoted $\,J_{x_0}^r \C \,$.

Let $\,\mbox{Diff}_{x_0}\,$ be the group of germs of local
diffeomorphisms of $\,X\,$ leaving $\,x_0\,$ fixed, and let
$\,\mbox{Diff}_{x_0}^r\,$ be the Lie group of $\,r-$jets at
$\,x_0\,$ of local diffeomorphisms of $\,X\,$ leaving $\,x_0\,$
fixed. We have the following exact sequence of groups:
$$0\longrightarrow H_{x_0}^r \longrightarrow \mbox{Diff}_{x_0}
\longrightarrow \mbox{Diff}_{x_0}^r \longrightarrow 0\,,$$
$H_{x_0}^r\,$ being the subgroup of $\,\mbox{Diff}_{x_0}\,$ made
up of those diffeomorphisms whose $\,r-$jet at $\,x_0\,$ coincides
with that of the identity.

The group $\,\mbox{Diff}_{x_0}\,$ acts on
$\,J_{x_0}^r \C$: if $\tau \in \mbox{Diff}_{x_0}$ and $\,j^r_{x_0}\nabla\in J^r_{x_0} \C$, then $\,\tau \cdot  (j^r_{x_0} \nabla)\,$ 
is the $r$-jet at $\,x_0\,$ of the linear connection $\,\tau \cdot \nabla$, defined as:

$$ (\tau \cdot \nabla)_D \bar{D} := \tau_*^{-1} \left( \nabla_{\tau_* D} (\tau_* \bar{D}) \right) \ . $$

Note that the subgroup $\,H_{x_0}^{r+2}\,$ acts
trivially on $\,J_{x_0}^r \C$, so the action of $\,\mbox{Diff}_{x_0}\,$ factors through an action of
$\,\mbox{Diff}_{x_0}^{r+2}\,$.

\smallskip

\begin{defi}
Two $\,r-$jets $\,j_{x_0}^r \nabla\,$,$\,j_{x_0}^r\bar{\nabla}\in
J_{x_0}^r \C\,$ are said to be \textbf{equivalent} if there exists a
local diffeomorphism $\,\tau\in \mbox{Diff}_{x_0}\,$ such that
$\,j_{x_0}^r\bar{\nabla}=j_{x_0}^r(\tau^*\nabla)\,$.
\end{defi}

\medskip

Equivalence classes of $\,r-$jets of linear connections constitute a ringed
space. To be precise:

\smallskip

\begin{defi}
We call \textbf{moduli space} of $\,r-$jets of linear connections the quotient 
ringed space
$$
\mathfrak{C}^r_n:=J_{x_0}^r \C/\mbox{Diff}_{x_0} =
J_{x_0}^r \C/\mbox{Diff}_{x_0}^{r+2}\,.
$$

In the case of symmetric connections,
the moduli space will be denoted $\,\widetilde{\mathfrak{C}}_{n}^r\,$.
\end{defi}

\medskip

The moduli space depends neither
on the point $\,x_0\,$ nor on the chosen $\,n-$dimensional
manifold.


\smallskip

\begin{ejem} \label{moduliendimensionbaja} If $\,n=1\,$, any linear connection is locally isomorphic to the standard, flat connection on $\,\RR$. Hence, all moduli spaces $\,\mathfrak{C}^r_1\,$ reduce to a single point, for any $\,r \in \mathbb{N} \cup \{0\}$.

If $\,\nabla\,$ is a symmetric connection around $\,x_0$, there always exists a chart in which $\, \Gamma_{ij}^k (x_0) = 0$ (see Corollary \ref{SimetriaEnElPunto}). Therefore, the moduli space of 0-jets of symmetric connections $\widetilde{\mathfrak{C}}^0_n\,$ also reduces to a single point, on any dimension.
\end{ejem}

\section{Description via normal tensors}

\begin{defi}
Let $\,m \geq 0\,$ be a fixed integer and let $\,x \in X\,$ be a point. The space of
\textbf{normal tensors} of order $\,m\,$ at $\,x\,$, which we
will denote by $\,C_m\,$, is the vector space of
$\,(1,m+2)-$tensors $\,T\,$ at $\,x\,$ having the
following symmetries:

- they are symmetric in the last $\,m\,$ covariant indices:
\begin{equation}\label{S1}T^l_{ijk_1
\ldots k_m} = T^l_{ijk_{\sigma(1)} \ldots k_{\sigma(m)}}\quad, \quad
\forall \ \sigma \in S_m\,;
\end{equation}


- the symmetrization over the $\,m+2\,$ covariant indices is zero:
\begin{equation}\label{S2}
\sum_{\sigma \in S_{m+2}} T^l_{\sigma (i) \sigma (j) \sigma (k_1)  \ldots \sigma (k_m)} = 0\,.
\end{equation}

When dealing with symmetric connections, we will consider a slightly different definition of normal tensors of order $m\,$ at $x\,$. The corresponding vector space will be denoted by $\widetilde{C}_m$, and it will consist of all $(1,m+2)-$tensors $T\,$ at $x\,$ verifying symmetries \ref{S1}, \ref{S2} and that of being symmetric in the first two covariant indices:
\begin{equation}
T_{ijk_1\ldots k_m}^l=T_{jik_1\ldots k_m}^l\,.
\end{equation}

Due to this additional symmetry, it is easily checked that $\, \widetilde{C}_0 = 0$.
\end{defi}

\bigskip

To show how a germ of linear connection $\, \nabla \,$ around $\,x\,$ produces a sequence
of normal tensors $\, \Gamma^m \,$ at $\,x\,$, let us briefly recall some definitions and results.

\begin{defi}
A chart $\,(x_1, \ldots , x_n)\,$ in a 
neighbourhood of $\,x\,$ is said to be a \textbf{normal} system for $\,\nabla \,$ at the point $\,x\,$ if the
geodesics passing through $\,x\,$ at $\,t=0\,$ are precisely the
``straight lines'' $\,\{ x_1(t) = \lambda_1 t , \ldots , x_n(t) =
\lambda_n t \}\,$, where $\,\lambda_i \in \mathbb{R}\,$.
\end{defi}

\medskip

As it is well known, via the exponential map
$\,\mbox{exp}_\nabla :T_{x}X\rightarrow X\,$, normal
systems on $\,X\,$ correspond bijectively to linear
systems on $\,T_{x}X\,$. Therefore, two normal systems differ in
a linear  transformation.

A simple, standard calculation allows to prove:

\begin{lema}
Let $\,(x_1,\ldots,x_n)\,$ be germs of a chart centred at $\,x\in X\,$, and let $\Gamma^k_{ij}$ be the Christoffel symbols of a linear conection $\nabla$ in those coordinates. 

It holds:
$$ (x_1, \ldots , x_n) \mbox{ is a normal system for } \nabla \quad \Leftrightarrow \quad \sum_{i,j=1}^n x_i x_j \Gamma_{ij}^k \, = \, 0 \quad , \ \  \, k = 1 , \ldots , n. \ . $$
\end{lema}

\medskip
Recall the exponential map $\mathrm{exp}_{\nabla} \colon T_x X \to X$ is a diffeomorphism around the origin. Let $\overline{\nabla}$ be the germ of linear connection around $x$ that corresponds, via the exponential map, to the canonical flat connection of $T_xX$. 

Let us also consider the difference tensor between $\nabla$ and $\overline{\nabla}$:
$$ \TT (\omega, D_1, D_2 ) := \omega \left( D_1^\nabla D_2 - D_1 ^{\overline{\nabla}} D_2 \right) \ . $$

If $(x_1, \ldots , x_n)$ is a normal chart for $\nabla$ around $x$, then:
$$ \TT := \sum_{i,j,k} \ \Gamma_{ij}^k \ \frac{\partial }{\partial x_k} \otimes \d x^i \otimes \d x^j  \ .  $$


\begin{defi} For each $m \geq 0$, the $\,m-$th \textbf{normal tensor} of the connection $\,\nabla \,$ at the point $\,x\,$ is:
$$ \Gamma^m_x \, := \, \overline{\nabla}^m_x \TT \ . $$

If $(x_1, \ldots , x_n)$ is a normal chart for $\nabla$ around $x$, then:
$$ \Gamma^m_x \, = \, \sum_{i,j,k , a_1, \ldots a_m} \Gamma_{ij; a_1 \ldots a_m}^k (0) \ \frac{\partial }{\partial x_k} \otimes \d x^i \otimes \d x^j  \otimes \d x^{a_1} \otimes \ldots \otimes \d x^{a_m} \ , $$
where 
$$ \Gamma_{ij; a_1 \ldots a_m}^k  := \frac{\partial \Gamma_{ij}^k }{\partial x_{a_1} \ldots \partial x_{a_m}}  \ . $$
\end{defi}

\smallskip
\begin{propo}
For each $\, m\geq 0$, the tensor $\,\Gamma^m_x\,$ belongs to $\,C_m$.
\end{propo}

\demo We only have to check that the symmetrization of the $m+2$-covariant indices of $\Gamma^m_x$ is zero. To this end, let $(x_1, \ldots , x_n)$ be a normal chart for $\nabla$ around $x$, so that the Christoffel symbols of $\nabla$ in these coordinates satisfy:

\begin{equation*}
 \sum_{i,j=1}^n x_i x_j \Gamma_{ij}^k \, = \, 0  \ . 
\end{equation*}


If we differentiate $m+2$ times this equality and evaluate at the origin, it follows that, for any $a_1, \ldots a_{m+2} \in \{  1, \ldots , n \}$:
\begin{equation}\label{EnelPunto}
\sum \  \Gamma_{a_i a_j; a_1 \ldots a_{m+2} }^k (0) + \Gamma_{a_j a_i ; a_1 \ldots a_{m+2} }^k (0) \ = \ 0 
\end{equation}
where the sum is over all the possible $i \leq j$, taken among $\{ 1 , \ldots , m+2 \}$. As the functions $\Gamma_{ij; l_1 \ldots l_m}^k$ are symmetric in the last $\,m\,$ indices, the thesis follows.

\qed

Using equality (\ref{EnelPunto}), it immediately follows:

\begin{coro}\label{SimetriaEnElPunto}
If $\nabla$ is symmetric, then $\Gamma^0_x = 0$.
\end{coro}

\subsubsection*{Reduction theorem}

If $\, \nabla\,$ is a germ of linear connection around $\, x$, let us denote $\, (\Gamma^0_x , \ldots ,  \Gamma^m_x , \ldots )$ the sequence of its normal tensors at the point $\, x$. 
Observe that $\Gamma^m_x$ only depends on $\, j^m_x \nabla$.

\begin{lema}[Orbit reduction] Let $\,G\,$ be a Lie group acting on a smooth manifold $\,X$, and let $\, f \colon X \to Y\,$ be a surjective regular submersion.

If the orbits of $\,G\,$ are precisely the fibres of $\,f$, then the quotient $\,X/G\,$ is a smooth manifold, and the map $\,[x] \, \longmapsto \, f(x)\,$ establishes an isomorphism of smooth manifolds:
$$ X/G \xrightarrow{\ \sim \ } Y \ . $$
\end{lema}

\proof The universal quotient property assures 
that $\,f\,$ factors through a unique morphism of ringed spaces $\,\bar{f} \colon X/G \to Y$, $\, [x] \to f(x)$, that is clearly bijective. Moreover, 
any local section $\,s\,$ of $\,f$ induces, when projected to the quotient, a morphism of ringed spaces that is a (local) inverse for $\,\bar{f}$. That is to say,
the map $\, \bar{f}\,$ is in fact an isomorphism of smooth manifolds.

\qed

\begin{teo}[Reduction]\label{LemaFundamental}
For each $\, r \in \NN \cup \{ 0 \}$, the map 
$$ J_x^r \Conn \ \xrightarrow{\ \ \pi_r \  \ } \ C_0 \times \stackrel{}{\ldots} \times C_r  \qquad , \qquad j^r_x \nabla  \ \to \ (\Gamma^0_x , \ldots ,  \Gamma^r_x )  $$ 
is a surjective regular submersion, whose fibers are the orbits of $H^{r+2}_x$.

Therefore, $\pi_r$ induces an isomorphism of smooth manifolds:
$$
\begin{CD}
(J^r_x \Conn) \, / \, H^{r+2}_x \ @= \ C_0 \times \ldots \times C_r \ .
\end{CD} $$
\end{teo}

\demo To check it is a regular submersion, let us construct a global section passing through any point of $J^r_x \Conn$. Such a section depends on a choice
of coordinates, and its image will be those jets having the chosen coordinates as a normal chart.

So let $\, (x_1, \ldots , x_n)\, $ be coordinates around $\, x$. Let us define a map 
$$\,s_r \colon C_0 \times \ldots \times C_r \longrightarrow J^r_x \Conn \quad , \quad s_r (T^0, \ldots , T^r) := j^r_x \nabla \ , $$ 
where $\, \nabla\,$ is the linear connection whose Christoffel symbols on the chosen coordinates are the following polynomial functions:
$$ \Gamma_{ij}^k \, := \,  (T^0)_{ij}^k \, + \, \sum_{a_1} (T^1)^k_{ij,a_1} x_{a_1} \, + \, \ldots \, + \, \frac{1}{r!} \sum_{a_1 \ldots a_r} (T^r)^k_{ij, a_1 \ldots a_r} x_{a_1} \ldots x_{a_r} \ .  $$

This map $\,s_r\,$ is clearly smooth and satisfies:

\begin{itemize}
\item It is a section of $\pi_r$: the chart $\, (x_1 ,\ldots , x_n)\, $ is a normal system for $\, \nabla\, $ around $\, x\,$, because the functions $\,x^ix^j \Gamma_{ij}^k\,$ vanish, due to the symmetries of the $\,T^m$. Therefore, the $\,r\,$ first normal tensors associated to $\,j^r_x \nabla\,$ at the point $\, x\,$ are precisely $\,T^0, \ldots , T^r$. 

\item This section can pass through any point $\, j^r_x \overline{\nabla}$, by simply choosing $\, (x_1, \ldots , x_n)\,$ to be a normal system for $\, \bar{\nabla}$.
\end{itemize}

\medskip Let us now check that the fibres of $\, \pi_r\,$ are the orbits of $\, H^{r+2}_x$. On the one hand, normal tensors are natural (i.e., independent of choices of coordinates), so that $\, \pi_r\,$ is $\,\Diff^{r+2}$-equivariant. Hence, as $\, H^{r+2}_x\,$ acts trivially on the spaces of normal tensors $\,C_0 \times \ldots \times C_r$, the orbits of $\, H^{r+2}_x\,$ are inside the fibres of $\,\pi_r$.

On the other hand, let $\,j^r_x \nabla\,$ and $\,j^r_x \overline{\nabla}\,$ be two jets of linear connections with the same normal tensors $\, \Gamma^0, \ldots , \Gamma^r\,$ at the point $\,x$. 

Let us fix a basis of $\, T_xX\,$ and let $\,(x_1, \ldots , x_n)\,$ and $\, (\overline{x}_1, \ldots \overline{x}_n)\,$ be the corresponding normal systems induced by those jets. 

Let $\, \tau\,$ the diffeomorphism carrying one chart to the other, $\,\tau(x_i) := \overline{x}_i$. As $\, \d_x x_i = \d_x \overline{x}_i$, because both coincide with the chosen basis,  it follows that $\, j^{r+2}_x\tau \in H_x^{r+2}$.
 
Now, an easy computation in coordinates allows to conclude that $\, \tau_* \left( j^r_x \nabla \right) = j^r_x \overline{\nabla}$, so that both jets are in the same orbit of $\, H^{r+2}_x$.

\qed

\begin{obser}
A similar argument proves that, for each $\, r \in \NN \cup \{ 0 \}$, the map 
$$ J_x^r \widetilde{\Conn} \ \xrightarrow{\ \ \pi_r \  \ } \ \widetilde{C}_0 \times \stackrel{}{\ldots} \times \widetilde{C}_r  \qquad , \qquad j^r_x \nabla  \ \to \ (\Gamma^0_x , \ldots ,  \Gamma^r_x )  $$ 
is a surjective regular submersion, whose fibers are the orbits of $H^{r+2}_x$.

Therefore, $\pi_r$ induces an isomorphism of smooth manifolds:
$$
\begin{CD}
(J^r_x \widetilde{\Conn}) \, / \, H^{r+2}_x \ @= \ \widetilde{C}_0 \times \ldots \times \widetilde{C}_r \ .
\end{CD} $$
\end{obser}

\begin{obser}
Fix a chart $\, (x_1, \ldots , x_n)\, $ be a chart around $\, x$, and let $\, \mathcal{N}^r_x\subset J^r_x \mathcal{C}\,$ be the submanifold formed by those jets for which $\, (x_1, \ldots , x_n)\, $ are normal coordinates.

The proof of the previous Theorem also says that $\,\mathcal{N}_x\,$ is a slice of the action of $\,H^{r+2}_x\,$ on $\,J^r_x \C$.
\end{obser}

\smallskip

Taking into account the exact sequence:

$$1 \longrightarrow H_{x}^{r+2} \longrightarrow \mbox{Diff}^{r+2}_{x}
\longrightarrow \Gl \longrightarrow 1\,,$$
the Reduction Theorem \ref{LemaFundamental} has the following immediate consequence, which is the main result of this paper:

\begin{teo}\label{elmoduli}
The moduli space of jets of linear connections is isomorphic, as a ringed space, to the orbit space of a linear representation of $\,\Gl $:
$$ \mathfrak{C}^r_n \, \simeq \, \left( C_0 \times \ldots \times C_r \right) / \Gl \qquad , \qquad \widetilde{\mathfrak{C}}^r_n \, \simeq \, \left( \widetilde{C}_0 \times \ldots \times \widetilde{C}_r \right) / \Gl  . $$
\end{teo}

Compare this result with similar statements obtained for Riemannian metrics (\cite{GNS}) and other $\, G-$structures (\cite{Jaime}). 



\section{Some properties of the moduli spaces}

This last section is devoted to extract some consequences of Theorem \ref{elmoduli}.

\subsection{Non-existence of differential invariants}

Let us consider the quotient morphism
$$
\xymatrix{J_{x}^r \C \ar[r]^-{\pi} & J_{x}^r \C/\mbox{Diff}_{x}\,=\,\mathfrak{C}^r_n\,.}
$$

\smallskip

\begin{defi}\label{definiciondeinvariante}
A (scalar) \textbf{differential invariant} of order $\,\leq r \,$ of
linear connections is defined to be
a global differentiable function on some $\,\mathfrak{C}^r_n\,$.
\end{defi}

\medskip

Taking into account the ringed space structure of
$\,\mathfrak{C}^r_n\,$ (see Example \ref{ExQuotient}), we can simply write:
$$
\{\mbox{Differential invariants of order}\,\leq r\}=\mathcal{C}^{\infty}(\mathfrak{C}^r_n)=\mathcal{C}^{\infty}(J_{x}^r \C)^{\mbox{Diff}_{x}}\,.
$$

\smallskip

\begin{lema}
For all $\,r \in \mathbb{N} \cup \{0\}$, the algebra of $\,\Gl$-invariant, polynomial functions 
$$ C_0 \times \ldots \times C_r \, \longrightarrow \, \RR $$ is trivial; i.e., it consists on constant functions only.
\end{lema}

\demo If a polynomial function is $\, \Gl$-invariant, then so they are its homogeneous components; hence, it is enough to argue the case of homogeneous polynomials.

The vector space of $\,\Gl$-invariant, polynomial functions $\,C_0 \times \ldots \times C_r \, \longrightarrow \, \RR$, homogeneous of degree $\,k\,$ is isomorphic to:
$$\bigoplus_{d_0 + \ldots + d_r = k} \mathrm{Hom}_{\Gl} \left( S^{d_0} C_0 \otimes \ldots \otimes S^{d_r} C_r \, , \, \RR \right) \ . $$

By Proposition \ref{Semisimple}, any $\,\Gl$-invariant linear map $\,S^{d_0} C_0 \otimes \ldots \otimes S^{d_r} C_r  \to \RR\,$ is the restriction of a $\,\Gl$-invariant linear map 
$$ \otimes T^*_xX \otimes \stackrel{p}{\ldots} \otimes T^*_xX \otimes T_xX \otimes \stackrel{q}{\ldots} \otimes T_xX \, \longrightarrow \, \RR \ ,$$ where $\,p = 2 d_0 + \ldots + (r+2) d_r$, and $\, q = d_0 + \ldots + d_r $.

If $\, k > 0$, then $\, p \neq q\,$ and Theorem \ref{MainTheoremGl} says that there are no such linear maps. That is to say, if $\, k > 0\, $ the above vector spaces reduce to zero and the thesis follows. 

\qed

\begin{teo} {\bf(Non-existence of differential invariants)}\label{finitud} The only differential invariants associated 
to (symmetric or not) linear connections  are constant functions.

That is to say,
$$ \mathcal{C}^\infty (\mathfrak{C}^r_n ) \, = \, \RR \quad , \quad \mathcal{C}^\infty (\widetilde{\mathfrak{C}}^r_n ) \, = \, \RR \ . $$ 
\end{teo}

\medskip

\proof By Corollary \ref{elmoduli} and the universal property of quotient ringed spaces,
$$
\mathcal{C}^{\infty}(\mathfrak{C}^r_n)\,=\, \mathcal{C}^\infty \left( ( C_0 \times \ldots \times C_r ) / \Gl \right) \, = \, \mathcal{C}^{\infty}(C_0\times \ldots \times C_r)^{\Gl}\,.
$$

Luna's theorem \ref{TeoremaLuna} describes such an algebra in terms of a system of generators of the algebra of polynomial, $\, \Gl$-invariant functions $\, C_0 \times \ldots \times C_r \to \RR$.

But the previous Lemma proves that any such a polynomial function is constant, and hence the algebra under consideration is trivial; i. e. $\,\mathcal{C}^{\infty}(\mathfrak{C}^r_n)\, = \, \RR$. 

An analogous reasoning applies for the case of symmetric connections, and proves $\,\mathcal{C}^{\infty}(\widetilde{\mathfrak{C}}^r_n)\, = \, \RR$.

\qed

\medskip
\begin{obser} More generally, any tensorial invariant, not necessarily scalar, is called a {\it natural tensor} associated to linear connections.
These natural tensors are usually described in terms of the curvature operator and its covariant derivatives, see (\cite{KMS}).

A similar argument to that presented above allows to produce an alternative description, using normal tensors, of the vector space of $\,(p,q)$-natural tensors (of order $\,\leq r$) associated to linear connections:

$$
\left\{\begin{aligned}\text{Smooth, $\,\Diff$-equivariant maps }\\
T\colon J^r_x \mathcal{C} \, \longrightarrow \, \otimes^p T^*_xX \otimes^q T_xX
\quad\end{aligned}\right\}$$
$$\begin{CD}
 @| \\
 \bigoplus \limits_{d_i} \, \mathrm{Hom}_{\Gl} \left(  S^{d_0} C_0 \otimes \cdots \otimes
S^{d_r} C_r  \ , \  \otimes^p T^*_x X \otimes^q T_x X \  \right)
\end{CD} $$  
where the summation is over all sequences $\,\{d_0 , \ldots , d_r \}\, $ of non-neagtive integers satisfying:
\begin{equation}\label{CondicionHomogeneidad}
\ \ d_0 + 2 d_1 + \ldots + (r+1) \, d_r  = p - q  \ . \end{equation}


As an application, a simple reasoning using Lemma \ref{Semisimple} and Theorem \ref{MainTheoremGl} allows to prove the following characterization of the curvature tensor of symmetric, linear connections (see \cite{KMS}, Section 28 for related results):

\medskip
{\it Up to constant multiples, the curvature tensor $\,R\,$ is the only natural 2-form with values on $\,\mathrm{End} (TX)\,$ associated to symmetric, linear connections.}
\end{obser}

\smallskip

\subsection{A few comments on dimensions of generic strata}




Recall that, due to Theorem \ref{elmoduli}, the following isomorphism of ringed spaces holds:
$$
\widetilde{\mathfrak{C}}_n^r= (J_x^r \widetilde{C})/\mathrm{Diff}_x = (\widetilde{C}_1\times \ldots \times \widetilde{C}_r)/ \mathrm{Gl}_n\,.
$$

Let us make some comments on what could be called ``generic dimension'' of this orbit space. To be precise, we will check that the formula
\begin{equation}\label{FormulaAplicamos}
\sum_{m=1}^r \mathrm{dim}\, \widetilde{C}_m - (\mathrm{dim}\, \mathrm{Gl}_n -i)\,,
\end{equation}
where $i$ denotes the minimum dimension of the isotropy groups for the action of $\mathrm{Gl}_n\,$ on $\widetilde{C}_1\times \ldots \times \widetilde{C}_r\,$, recovers other formulae appearing in the literature regarding the dimension or the Poincar\'e series of the, loosely speaking, ``strata of generic jets'' (\cite{Arnold}, \cite{Dubrovskiy}).

If $s_{m+2}$ denotes the symmetrization operator, then the following sequence is exact:

\begin{equation*}\label{exactadesimetricas}
0 \longrightarrow \widetilde{C}_m \longrightarrow T_x X \otimes S^2 T^*_xX \otimes S^m T^*_xX \xrightarrow{\ \ s_{m+2} \ \ } T_xX \otimes S^{m+2} T^*_xX \longrightarrow 0 \,.
\end{equation*}

Using this sequence, a straightforward computation yields the dimension of $\widetilde{C}_m\,$:
$$ \mathrm{dim}\, \widetilde{C}_m = n \frac{n(n+1)}{2}{n+m-1 \choose m}-n{n+m+1 \choose m+2} \ . $$


\medskip

Later we will check that, if $\mathrm{dim}\, X=n= 2\,$,  then any 1-jet has, at least, a one-dimensional isotropy group; in any other case, generic jets have no isotropy. That is to say, we have $i=1\,$ if $(n,r) = (2,1)$, and $i=0\,$ otherwise. 

Therefore, taking into account that $\mathrm{dim}\,\mathrm{Gl}_n=n^2\,$, we observe that formula \ref{FormulaAplicamos} produces the same result about the generic dimension of the moduli space $\widetilde{\mathfrak{C}}_n^r\,$ that can be found in \cite{Dubrovskiy} for dimension $n\geq 2\,$ (the trivial case $n=1\,$ has already been dealt with in Example \ref{moduliendimensionbaja}):
\begin{align*}
\mathrm{dim}\, \widetilde{\mathfrak{C}}_n^r &=n \frac{n(n+1)}{2}\sum_{m=0}^{r}{n+m-1 \choose m}-n\sum_{m=0}^{r}{n+m+1 \choose m+2}- \left( n^2-\delta_2^n\delta_1^r \right) \\
&=n\frac{n(n+1)}{2}\sum_{m=0}^{r}{n+m-1 \choose n-1}-n\sum_{m=1}^{r+2}{n+m-1 \choose n-1}+\delta_2^n\delta_1^r\,,
\end{align*}
and, hence, it also provides exactly the same expression of the Poincar\'{e} series that can be read in \cite{Dubrovskiy} (page 1055).

\subsubsection*{Isotropy groups of 1-jets in dimension 2}\label{isotropia}

Finally, let us make some comments regarding the isotropy groups of generic 1-jets in dimension 2 (compare with \cite{Dubrovskiy}, where a similar goal is achieved through direct computation).

\begin{defi}
The vector space of {\bf curvature-like tensors} is the subspace $\, \R \ \subseteq \  \Lambda^2 T^*_xX \otimes T^*_xX \otimes T_xX\, $ defined by the linear Bianchi identity:
\begin{equation}\label{BianchiGamma}
R_{ijl}^k + R_{lij}^k + R_{jli}^k = 0 \ . 
\end{equation} 
\end{defi}

These curvature-like tensors are closely related to normal tensors: it is not difficult to check that the linear map $\, R^k_{ijl} := \Gamma_{jli}^k - \Gamma_{ilj}^k$ establishes an isomorphism of $\, \Gl$-modules:
$$ \widetilde{C}_1 \ \simeq \ \R \ , $$ whose inverse is $\, \Gamma_{ijl}^k := \frac{1}{3} \left( R^k_{lij} + R^k_{lji} \right)$.

\smallskip

Let us fix some notations: the symmetrization, skew-symmetrization and {\it Ricci} maps will be denoted, respectively,
$$ s \colon \otimes^2 T^*_xX \to S^2T^*_xX \quad , \quad a \colon \otimes^2 T^*_xX \to \Lambda^2 T^*_xX \quad , \quad \R \ \xrightarrow{\, \rho \,} \ \otimes^2 T^*_xX \ , $$
where $\, \rho(R)_{ij} := \sum_{k=1}^n R_{ikj}^k$. 


\begin{lema}[\cite{Gilkey}, Lemma 4.4.1]\label{LemaGilkey}
If $\,X\,$ has dimension 2, then the Ricci map establishes an isomorphism of $\, {\rm Gl}_2$-modules:
$$ \R \ \xrightarrow{\ \rho_s \oplus \rho_a \ } \ S^2 (T^*_xX) \oplus \Lambda^2 (T^*_xX) \ , $$
where $\,\rho_s := s \circ \rho\,$ and $\, \rho_a := a \circ \rho$.
\end{lema}

\smallskip

This Lemma implies that, if $\,X\,$ has dimension 2, the isotropy group of any 1-jet $\, j^1_x \nabla\,$ under the action of $\, \Diff\,$ is, at least, 1-dimensional. 

In fact, due to the isomorphisms
$$ \widetilde{\mathfrak{C}}^1_n \, = \, \widetilde{C}_1 / {\rm Gl}_2 \, = \, \R / {\rm Gl}_2 \, = \, (S^2 T^*_xX  \oplus \Lambda^2 T^*_xX ) / {\rm Gl}_2 \ , $$
it is enough to check that any pair $\, (T_2 , \omega_2)\,$ of a symmetric 2-tensor and a 2-form on a 2-dimensional vector space has, at least, a 1-dimensional isotropy group under the action of $\, {\rm Gl}_2$.

If the metric $\,T_2\,$ is non-singular, then its automorphisms have determinant equal to 1, and hence preserve any 2-form $\,\omega_2$.
In this case, the isotropy group of the pair $\,(T_2, \omega_2)\,$ is isomorphic to $\,O(2)\,$ or $\,O(1,1)$, depending on the signature of $\,T_2$.

The other cases where $\,T_2\,$ is singular are easily analyzed in a similar manner, resulting in larger isotropy groups.

Analogous arguments, with the corresponding versions of Lemma \ref{LemaGilkey}, can be applied to check that, if $\, n > 2\,$ or $\, k > 1$, then the isotropy group of a ``generic'' jet of linear connection is trivial.

\end{document}